\documentstyle [11pt]{article}
\begin{document}
\font\tfont= cmbx10 scaled \magstep3
\font\afont= cmcsc10 scaled \magstep2
\centerline{\tfont A Representation of Stable Banach Spaces}
\bigskip
\centerline{\afont Fouad Chaatit }
\bigskip
\begin{abstract}

{We show that any separable stable Banach space
can be represented as a group of isometries on a separable reflexive
Banach space, which extends a result of S. Guerre and M. Levy. 
As a consequence, we can then represent homeomorphically
its
space of types.}
\end{abstract}

\newtheorem{thm}{Theorem }[chapter]
\newtheorem{lem}[thm]{Lemma} 
\newtheorem{cor}[thm]{Corollary}
\newtheorem{prop}[thm]{Proposition}          
\newtheorem{sublem}[thm]{Sublemma}
\newtheorem{ex}{Exercise}                             
\newenvironment{proof}{\medskip \par \noindent {\bf Proof:}\ }{\hfill $\Box$ 
		       \medskip \par}

\setcounter{section}{0}
\setcounter{thm}{0}
\setcounter{ex}{0} 
\setcounter{page}{1}

\newcommand{\bdfn}{\noindent{\bf Definition:\,} \begin{em}}
\newcommand{\edfn}{\end{em}}
\newcommand{\ntn}{\noindent{\bf Notations:\,\,}}

\newcommand{\rmk}{\noindent{\bf Remark\,\,}}

\newcommand{\del}{\mbox{$\delta$}}
\newcommand{\lftyn}{\mbox{${l_{\infty}^{n}}$}}
\newcommand{\norm}[1]{\mbox{$\|#1\|$}}             
\newcommand{\som}{\mbox{$\displaystyle\sum_{n=0}^{\infty}$}}
\newcommand{\son}{\mbox{$\displaystyle\sum_{n=1}^{\infty}$}}
\newcommand{\reals}{\mbox{\bf R}}
\newcommand{\complexes}{\mbox{\bf C}}
\newcommand{\integers}{\mbox{\bf N}}
\newcommand{\smint}{\mbox{\scriptsize {\bf N}}}
\newcommand{\complex}{\mbox{\bf C}}
\newcommand{\ep}{\mbox{$\varepsilon$}}
\newcommand{\B}{\mbox{${\cal B}_1(K)$}}
\newcommand{\pf}{\noindent{\bf Proof:\,\ }}
\newcommand{\ra}{\mbox{$\longrightarrow$}}
\def\qed{~\hfill~\blackbox\medskip}
\def\blackbox{\hbox{\vrule width6pt height7pt depth1pt}}

\newcommand{\osm}{\mbox{$(\Omega,\Sigma,\mu)$}}
\newcommand{\om}{\mbox{$\omega$}}

\section{Introduction}

Stable Banach spaces have been introduced for the first time
 by J.L. Krivine and B. Maurey \cite{bib:K.M}.
In this  famous paper, a rather  simple condition on the norm is given to
ensure, for any $\ep > 0$ and any
infinite dimensional subspace $Y$  of a given Banach space $X,$  the
existence of a further subspace $Z$ of $Y$  that is  $(1 + \ep)$
isomorphic to $l_p,$ for some  $ 1 \leq p < \infty;$  i.e. there exists an
isomorphism $T: l_p \ra Z$ so that $\|T\| \|T^{-1}\| \leq (1+\ep).$
This generalized the result of D. Aldous which was  established when
$X= L_1$ and $1 \leq p \leq 2.$  Subsequently, this class has 
been intensively studied by various authors. We refer the reader 
to the book by S. Guerre \cite{bib:G} for an extensive  survey 
on stability in Banach spaces.
In this note, we are interested in getting a representation for 
separable stable Banach spaces. In \cite{bib:G.L}
S. Guerre and M. Levy characterized the real number $p$ so that a subspace $E$
of $L_1$ contains $l_p,$ as being the upper bound of the set of all  reals
$q$
so that $E$ embeds into $L_q.$ 
Actually, since this number $p(E)$ is also the upper bound of the set
of all reals $q$ so that $E$ is of  type $q$-Rademacher,
and B. Maurey and G. Pisier \cite{bib:M.P} have proved
that $l_{p(E)}$ is finitely representable in $E,$ 
S. Guerre and M. Levy prove that for any infinite dimensional
closed subspace $E$ of $L_1,$ if $l_{p}$ is finitely representable in
$E,$ then there exist a $q \leq p$ so that $E$ contains a subspace
isomorphic to $l_q.$
To do so, S. Guerre and M. Levy used a representation of the norm of
subspaces of $L_1$ and then a representation of the types on these
subspaces, (and on their ultraproducts)
as an inner product on some Hilbert space. Our results were inspired
from
this proof but they don't generalize to give a similar result
for a general stable Banach space. To see this one can
take for example the $l_p$-sum of $l_2^n$ which is stable, 
has $l_2^n$ finitely represented, but there exist no $q < p$
such that this space contains an isomorphic copy of $l_q.$
However we do obtain some generalizations of the techniques in
\cite{bib:G.L}.

 In Theorem 1, we obtain that any stable
separable Banach space $X$ can be represented homeomorphically
as a closed subgroup $G$, of the group of all isometries on a reflexive
separable Banach space $R,$ when $G$ is endowed with the strong
operator
 topology. This is used to represent the types on $X$ in the following
 way:
there exists a homeomorphism $T$ from ${F}_{X},$ the
Alexandroff compactification of ${\cal F}_X,$
the space of types on $X,$ endowed with topology of pointwise
convergence on $X,$ to $\tilde{G},$ closure for the weak operator
topology of $G$ so that
$T_{\sigma \star \tau} =  T_{\tau \star \sigma } = T_{\sigma} \circ
T_{\tau}$
for all types $\sigma$ and $\tau$ on $X.$
This will be Theorem .4. We also obtain that
$T$ is actually a homeomorphism between $({F}_{X}, s)$ and
$(\tilde{G}, \tau_{s})$ where $s$ is the topology of uniform
convergence on bounded subsets of $X$ and $\tau_{s}$ is the topology
of pointwise convergence on the reflexive Banach space $R.$ This will
be Theorem .7.
 
It would be interesting to investigate if there is a similar representation
of the norm of some non stable Banach 
spaces,  in particular when they
{\it don't} contain $l_{p}$ (Tsirelson's space, Schlumprecht's space
or Gower-Maurey's space), and on what kind of spaces do they have
such a
representation, if there is any.
\bigskip

\noindent{\bf Acknowledgement:} I would like to warmly thank Professor
Bernard Maurey for many fruitful discussions and suggestions regarding this
work.

\section{Definitions and Notations:} 
Throughout this paper the unit ball  $\{ x \in X: \|x\| \leq 1\}$ of a
Banach space $X$ will
be denoted by $Ba X.$ Let $X$ be a separable Banach space. A {\em type}
$\tau$ on $X$ is
a function $\tau: X \ra [0,\infty)$ which is realized by a sequence
$(x_n) \subset X,$ meaning that 
$$\tau (x) = \lim_{n \rightarrow \infty} \|x + x_n \|. \,\,(*)$$
The sequence $(x_n)$ is then called {\em type determining.}
Note that since  the space $X$ was taken to be {\em separable}
the original definition in \cite{bib:K.M}  and this one are the same
for anytime $x_n$ is a bounded sequence in $X$ there exists
a subsequence $y_n$ that satisfies (*). Indeed given $x_n$ bounded in $X$
and $D$ countable dense in $X$ choose $y_n$ a subsequence of $x_n$ so that
$\lim_{n \rightarrow \infty} \|x_n + d\|$ exists for all $d \in D,$ by a
diagonal argument.
Then $y_n$ generates a unique type $\tau.$

A Banach space $X$ is called stable if whenever $(x_n)$ and $(y_n)$
are type determining then
$$\lim_{n \rightarrow \infty} \lim_{m \rightarrow \infty} \|x_n + x_m\|
= \lim_{m \rightarrow \infty} \lim_{n \rightarrow \infty} \|x_n + x_m\|.$$
%These definitions and concepts were introduced in \cite{bib:K.M.}.
Let ${\cal F}_X$ be the set of all types on $X.$
Several topologies can be defined on ${\cal F}_X.$

{\em The weak topology w}: this topology is induced by the metric
$$d_w(\tau, \sigma) = \sum_{k=1}^{\infty} 2^{-k}\frac{ |\tau(r_k) -
\sigma(r_k)|}{\|r_k\|}$$
where $\{r_k\}$ is  a dense subset of $X.$
Note that $\tau_n$ converges to $\tau$ for the metric $d_w$
if and only if $\tau_n$ converges to $\tau$ pointwise on $X.$
Also, if $X$ is separable, the space ${\cal F}_X$ is locally compact
and separable \cite{bib:K.M}.

{\em The strong topology}: which is defined by 
$$d_s(\tau, \sigma) = \sum_{k=1}^{\infty}2^{-k} \sup_{\|x\| \leq
k}|\tau(x) - \sigma(x)|.$$
If $X$ is stable separable then  $({\cal F}_X, d_s)$ is separable, see
\cite{bib:O} or \cite{bib:R}.

{\em The uniform topology}: which is defined by 
$$d_u(\tau, \sigma) = \sup_{x \in X} |\tau(x) - \sigma(x)|.$$
If $X$ is uniformly stable (see \cite{bib:C}) then $({\cal F}_X, d_u)$ is
separable, 
however ${\cal F}_{L_1}$ \cite{bib:H}, the space of types on $L_1$, 
 as well as  ${\cal F}_{Ts},$\cite{bib:O}
the space of types on Tsirelson's space,  are separable. 

A degenerate type $\tau_x$ is a map from $X$ to $[0,\infty)$ defined
by
$\tau_x(y) = \|x+y\|.$ Note that the space of all degenerate types is
naturally
homeomorphic to $X$ when endowed with the pointwise convergence
topology.
Usually the space of types on $X$ is viewed as the closure of the set
$\{\tau_x: x \in X\}$ in the product topology of $R^{X}.$
In case $X$ is separable, for all type $\tau$ on $X$, there exists a
sequence $(x_n)$ in $X$ so that $\|x_n + y\|$ converges to $\tau (y)$
for all $y \in X.$

For $\sigma$ and $\tau$ two types on a stable separable Banach space  $X$
the  convolution
is defined by:
$$\sigma \star \tau (x) = \lim_{n \rightarrow \infty} \lim_{m
\rightarrow \infty} \|x + x_n + y_m \|$$
where $(x_n)$ determines $\sigma$ and $(y_m)$ determines $\tau$.
We would like to point out that the hypothesis of stability of $X$
is essential here, and that is the reason why these results don't extend to 
general Banach spaces.

Under the stability hypothesis, the bracket of $\sigma$ and $\tau$ is
defined by:
$$[\sigma,\tau] = \sigma \star \tau (0)$$

and the norm of $\sigma$ is defined by:
$$\|\sigma\| = \sigma (0).$$
Let us now recall some facts about Lions-Peetre interpolation spaces. We
refer the reader to  \cite{bib:L.P} or \cite{bib:B} for an exhaustive study
of these spaces:

An ordered pair $(A_0,A_1)$ of Banach spaces is an {\it interpolation pair}
if these Banach spaces are given as subspaces of a common Hausdorff
topological linear space $Z$ so that the embeddings of $A_0$ and $A_1$
in $Z$ are continuous. In particular if $A_0$ embeds in $A_1$ then one can
take $A_1$ as the common space.

Given an interpolation pair $(A_0,A_1), \,\,0 <\theta<1\,\,$ and $1 \leq p
\leq \infty$, the interpolation space $[A_0,A_1]_{\theta,p}$ is the space of
all $z \in A_0 + A_1$ which admit  a representation $z=x_1(t) +x_2(t),
\,\,t\in \reals$ so that $e^{\theta t} x_1(t) \in L_p(A_0)$ and
$e^{-(1-\theta)t}x_2(t) \in L_p(A_1).$
The norm on this space is defined by
$$\|z\|_{\theta,p} = \inf \{{\|e^{\theta t} x_1(t)\|}^{1-\theta}_{p}
{\|e^{-(1-\theta)t} x_2(t)\|}^{\theta}_{p} \,\,\mbox{over all
possible}\,\,z= x_1(t)+x_2(t) \,\,t \in \reals\}$$
This norm satisfies  the following important inequality which we will refer
to as the interpolation inequality:
For all $x\in A_0 \cap A_1$ we have 
\begin{equation} 
\|x\|_{\theta,p} \leq C(\theta,p) \|x\|^{1-\theta}_{A_0} \|x\|^{\theta}_{A_1}
\end{equation}
\noindent{where} $C(\theta,p)$ is an absolute constant depending only on
$\theta$ and $p$. This follows from the fact that $z= \psi(t)z + (1-\psi(t))z$
where $\psi (t)$ is so that $0 \leq \psi(t) \leq 1,\,\, e^{\theta t}\psi(t)
\in L_p(\reals)$ and  $e^{-(1-\theta)t} \psi(t) \in L_p(\reals).$

An operator $T$ is said to be {\it bounded from the interpolation pair
$(A_0,A_1)$ to the interpolation pair $(B_0,B_1)$} if $T$ is defined on the
sum $A_0 +A_1$ 
and is a bounded operator from $A_0$ to $B_0$ and from $A_1$ to $B_1.$

Let us now recall the important interpolation Theorem (see \cite{bib:L.P} or
\cite{bib:B}):
\bigskip

\noindent{\bf Interpolation Theorem:} 
{\it Let $T$ be a bounded operator from the interpolation pair $A_0,A_1)$ to
the interpolation pair $(B_0,B_1).$ Then, for any $0<\theta<1$ and  $1 \leq
p \leq \infty$  we have that $T$ is  a bounded operator from
$[A_0,A_1]_{\theta,p}$
to  $[B_0,B_1]_{\theta,p};$ furthermore, $$\|T\|_{\theta,p}
\leq\|T\|^{1-\theta}_{0} \|T\|^{\theta}_{1}$$ where $\|T\|_{i}\,\,i=0,1$ is
the norm of $T$
as an operator from $A_i$ to $B_i.$}
\bigskip

\noindent{and}  the density theorem: \cite{bib:B}
\bigskip

\noindent{\bf Density Theorem:} 
{If $0< \theta <1$ and $1\leq p \leq \infty$ then $A_0 \cap A_1$ is dense in
$[A_0,A_1]_{\theta,p}.$} 
\bigskip

\section{Main result:}
\begin{thm}
Let $X$ be a separable stable Banach space.
There exists a representation $T$ from $X$ onto a closed subgroup
$G$ of  the group ${\cal U}$ of all the isometries on a separable reflexive
space $R,$
where ${\cal U}$ is endowed with the strong operator topology, so that
$T_{x+y} = T_x \circ T_y$ and $(T_x)^{-1} =T_{-x}$.
\end{thm}

Here the word representation means a group isomorphism that is bicontinuous
when 
$X$ is endowed with the strong operator topology and ${\cal U}$ is endowed
with the strong operator topology.

To prove this theorem, we first need to establish the following
proposition:

\begin{prop}
Let $X$ be a separable stable Banach space, and define
an operator $\Phi: {\cal M}(X) \ra {\cal C}_b(X) $
from the space of Radon measures on $X$ to the space of bounded continuous
functions
on $X$ by: $$\Phi(\mu)(y) = \int_{X} e^{-\|x-y\|} \mu(dx)$$
for all $\mu \in {\cal M}(X)$ and $y \in X.$
Then the operator $\Phi$ is weakly compact.
\end{prop}

 We first establish  the following

\begin{lem}
Let $F_X$ be the Alexandroff compactification 
of ${\cal F}_X,$ the space of all types on $X;$
and let $\widehat{f}: { F}_X  \times {F}_X  \ra \reals$ defined by
$\widehat{f}(\sigma, \tau)= \exp (-[\sigma,-\tau])$ if both $\sigma$ and
$\tau$
are in ${\cal F}_X,$ and $\widehat{f}(\sigma, \tau) =0$ otherwise.
Then $\widehat{f}$ is separately continuous.
\end{lem}

\pf 
Let $f:X \times X \ra \reals$ defined by: $f(x,y) = e^{-\|x-y\|}$ for $x, y$
in $X.$
It is clear that $f$ is then separately continuous and satisfies
the inequalities
$$f(x,y) \leq e^{-\|x\|} e^{\|y\|} \,\,\mbox{and}\,\, f(x,y) \leq e^{\|x\|}
e^{-\|y\|}.$$
By density of $X$ in ${\cal F}_X,$ we have that
$\widehat{f}: {{\cal F}_X  \times{\cal F}_X } \ra \reals$ defined by
$\widehat{f}(\sigma, \tau)= \exp (-[\sigma,-\tau])$ is separately continuous,
since the bracket [.,.] is, and
satisfies the inequalities:
$$\widehat{f}(\sigma, \tau) \leq \exp (-\|\sigma\|) \exp\|{\tau}\|
\,\,\mbox{and}
\,\,\widehat{f}(\sigma, \tau) \leq \exp {\|\sigma\|} \exp{(-\|\tau\|)}.$$
Using the last inequalities one can then extend $\widehat{f}$ to $F_X \times
F_Y$
by setting $\widehat{f}({\sigma}_{\infty},\tau) = 0 $ if ${\sigma}_{\infty}$
is the "point at infinity", $\widehat{f}({\sigma},{\tau}_{\infty}) =0$ if
${\tau}_{\infty}$
is the "point at infinity" and
$\widehat{f}({\sigma}_{\infty},{\tau}_{\infty}) = \lim_{\tau \ra
{\tau}_{\infty}} \widehat{f}({\sigma}_{\infty},\tau) = 0.$
\qed

We are now ready to prove Proposition .2. These ideas were already contained in 
\cite{bib:K.M} and \cite{bib:R}, however we include the  proof here for the
sake of completeness.

\noindent{\bf Proof of Proposition.2} Let $ \widehat{\Phi}$ be defined on $
{\cal M}(F_X)$
the space of Radon measures on $F_X,$ with values in ${\cal C} (F_X),$
the  space of continuous functions on $F_X,$ by $$ [\widehat{\Phi} (\mu)] (\tau)
= \int_{F_X} \widehat{f} (\sigma, \tau) \mu(d \sigma).$$
Then the operator $\Phi (\mu)$ is just the restriction of  $\widehat{\Phi}
(\mu)$ to
the degenerate types. 
On the other hand,  $\widehat{f}(\sigma, \tau)$ is a separately continuous map
on a product of two compact metrizable spaces, and thus is  of the first
Baire class on the product $F_X \times F_X$  \cite{bib:K.M}. It then follows
that $g$  is Borel on $F_X \times F_X,$  and by Fubini's theorem, if $\nu
\in {\cal C}(F_X)^{*}$ then
$\widehat{\Phi}^{*} (\nu) (\sigma) = \int_{F_X} \exp (-[\sigma, -\tau]) \nu
(d \tau)$
for $\sigma \in F_X.$ Using Lebesgue dominated convergence theorem,
we get that $\widehat{\Phi}^{*} (\nu)$ is actually in ${\cal C}(F_X)$
i.e  $\widehat{\Phi}^{*} :{\cal C}(F_X)^{*} \ra {\cal C}(F_X),$
an equivalent way of saying that  ${\widehat{\Phi}}^*$ is  a weakly compact
operator.
Now let $i$ be the canonical embedding of ${\cal M}(X),$  the space of Radon
measures on $X,$ in  $ {\cal M}(F_X) =  {\cal C}(F_X)^{*};$ and let $r$
be the operator from ${\cal C}(F_X)$ to the space of bounded continuous
functions ${\cal C}_{b}(X)$ on $X.$
Then $\Phi  = r \circ \widehat{\Phi} \circ i$ and this proves the
proposition.\qed

\noindent{\bf Proof of the Main Result:}
Let $\Phi : {\cal M}(X) \ra {\cal C}_{b}(X)$ be the operator defined in 
the previous proposition. By a result of W.J. Davis T. Figiel W.B. Johnson
and A. Pelczynski \cite{bib:D.F.J.P} the weakly compact operator $\Phi$
factors through a reflexive space $R.$ More precisely, B. Beauzamy proved in 
\cite{bib:B} that if a Banach space $A_0$  embeds in  a Banach space $A_1,$
then the interpolation spaces $[A_0,A_1]_{\theta,p}$ for any $0 < \theta <
1$ and $1<p<\infty$ are reflexive if and only if the injection of $A_0$ in $A_1$
is weakly compact.
Let $A_0 = \Phi ({\cal M}(X))$ endowed with the norm defined by the gauge of
the closed bounded convex symmetric set $\Phi (Ba {\cal M}(X)),$ and 
$A_1 = {\cal C}_{b} (X)$ with the uniform norm. Then $A_0$ is a Banach space 
whose unit ball is $\Phi (Ba {\cal M}(X)),$  is isometric to ${\cal
M}(X)/Ker \Phi,$ and the embedding of $A_0$ in $A_1$ is weakly compact.
The Banach spaces  $[A_0,A_1]_{\theta,p}$ for all $0 < \theta < 1$ and
$1<p<\infty$ are then  reflexive and $\Phi$ factors through these spaces since 
$\phi = j \circ i \circ \Phi$ where $i: A_{0} \ra [A_{0},A_{1}]_{\theta,p}$
and $j:[A_{0},A_{1}]_{\theta,p} \ra A_1$ are the canonical injections.

For each $z \in X,$ we define a translation operator on ${\cal C}_b (X)$
by $T_z f (x) = f(x-z)$ for all $x \in X.$ By duality we can also define an
operator on ${\cal M}(X),$
which we still denote  $T_z,$ by 
$$<T_z \mu , f> = <\mu , T_{-z} f> \,\, \mbox{for all}\,\, \mu \in {\cal M}(X)
\,\,\mbox{and}\,\,f \in {\cal C}_b (X).$$

We now verify  that $T_z$ commutes with $\Phi,$ i.e 
$T_z(\Phi (\mu)) = \Phi (T_z (\mu))$ for all $\mu \in {\cal M} (X).$
It suffices for that to see that
for all $x\in X$ and $y\in X$ we have
\begin{eqnarray*}
    T_{z} (\Phi({\delta}_{x})) (y)  & =  & \Phi({\delta}_{x}) (y-z) \\
				& =  & e^{-\|x -(y-z)\|} \\
				& =  & \Phi(\delta_{x+z} ) (y)\\
				& =  & \Phi (T_z(\delta_x))(y)
\end{eqnarray*}
The last equality  is easily checked since for all $f\in {\cal C}_b (X)$ we
have that 
$$<T_z\delta_x,f> = <\delta_x,T_-z f> = <\delta_x,f(.+z)> = <\delta_{x+z},f>.$$

Now the operator $T_z$ is an isometry of ${\cal C}_b (X),$ and when
restricted to  $A_0 = \Phi ({\cal M}(X))$ is an operator of norm $\|T_z\|_0
\leq 1.$ 
Using the Interpolation Theorem  stated in Section 2, $T_z$ is also an
operator on $[A_0,A_1]_{\theta,p}$ and its norm satisfies
$$\|T_z\|_{\theta,p} \leq \|T_z\|_{0}^{1-\theta} \|T_z\|_{1}^{\theta} \leq 1.$$
The same argument applied to the operator $T_{-z},$ inverse operator of $T_z,$
yields $\|T_{-z}\|_{\theta,p} \leq 1.$ In other words, $T_z$ is an isometry
of the reflexive space $R= [A_0,A_1]_{\theta,p}.$

Let now $G$ be the group of all the isometries $T_z,$ for $z \in X,$
endowed with  the strong operator topology $s.$ Suppose that $y_n$ converges
in norm to $y$ as $n \ra \infty$ and $r = \phi(\mu)$ is an arbitrary element
of $A_0.$ Since $T_y(\Phi (\mu)) = \Phi (T_y (\mu)),$ we have that $T_{y_{n}}r
- T_y(r) \in A_0$ and 
\begin{equation}
\|T_{y_{n}}r - T_y(r)\|_{A_0} \leq 2 \|r\|_{A_0}
\end{equation}
On the other hand, 
by definition of the norm in $A_1,$  
\begin{eqnarray*}
\|T_{y_{n}}r - T_{y}(r)\|_{A_1} & =  & \sup_{z \in X} |r(z-y_n) - r(z-y)| \\
				& =  & \sup_{z \in X}\{\phi (\mu + \nu) (z-y_n) - \Phi (\mu + \nu)(z-y)
\,\, \,\,:\nu \in Ker\Phi \}\\
			& \leq & \| \mu + \nu \| \, \| y_n - y\| \,\,\,\, 
			\mbox{for all}( \nu \in Ker\Phi 
\end{eqnarray*}

\noindent{But} since $A_0$ is isometric to ${\cal M}(X)/Ker \Phi,$
we have $\|\Phi (\mu)\|_{A_0} = \inf \{\|\mu + \nu \|, \Phi (\nu) = 0\},$
and so
\begin{equation}
\|T_{y_{n}}r - T_y(r)\|_{A_1} \leq \|r\|_{A_0} \|y_n - y\|
\end{equation}
Combining (2) ,(3) and the interpolation inequality  of section 2, we get 
\begin{equation}
\|T_{y_{n}}r - T_y(r)\|_{\theta, p} \leq C(\theta,p) 2^{1-\theta}
\|r\|_{A_0} \|y_n - y \|_{A_1}^{\theta}.
\end{equation}
Since the group $G$ is equi-continuous and $A_0$ is dense in
$[A_0,A_1]_{\theta,p},$ this inequality shows that
for any $r \in  [A_0,A_1]_{\theta,p}$ we have $T_{y_n}r \ra T_y r $ as $n\ra
\infty$
by Ascoli's theorem.

We now show that $G$ is sequentially closed.
Suppose indeed that $T_{y_n}$ converges strongly to an operator $S.$
For any $r \in R$ the sequence $T_{y_n} r$ is Cauchy, in particular when  $r
= \Phi (\delta_0).$ Using the embedding of $R$ in ${\cal C}_b (X)$
we get:
\begin{eqnarray*}
\|T_{y_n} r - T_{y_m} r\|_{R} & \geq & C^{-1} \|T_{y_n} r - T_{y_m}
r\|_{\infty} \\
& = & C^{-1} \left\|\exp{(-\|y_m-.\|)} - \exp{(-\|y_n -
.\|)}\right\|_{\infty} \\
&\geq & |\exp{(-\|y_n-y_m\|)} -1|
\end{eqnarray*}

\noindent{Thus} $y_n$ is Cauchy in $X;$ so letting $y = \lim_{n \ra \infty}
y_{n}$ and
using the previous arguments we get that $T_{y_n} r  \ra T_y r$ for any
$r \in R,$ and since also $T_{y_n} r  \ra S,$ we have that $T_y = S$ i.e
$G$ is sequentially closed.

We now show that the reflexive space $R = [A_0,A_1]_{\theta,p}$ is separable
provided $0 < \theta < 1$ and $1<p < \infty.$
Recall that $A_1 = {\cal C}_b(X)$ while $A_0 = \Phi({\cal M}(X)) = r \circ
\widehat{\Phi} \circ i ({\cal M}(X)) $ where
 $i$ is  the canonical embedding of ${\cal M}(X),$  the space of Radon
measures on $X,$ in  $ {\cal M}(F_X)$  and  $r$
 the  restriction to the types realized in $X$ defined from ${\cal C}(F_X)$
to the space of bounded continous functions ${\cal C}_{b}(X)$ on $X.$ Then
$i(A_0)$ is separable for the  ${\cal C}(F_X)$ norm, and therefore $A_0$ is
separable for  the ${\cal C}_{b}(X) = A_1$ norm. Let then $a_n$ be a dense
sequence in $A_0$ for the $A_1$ norm. 
Let $r \in Ba R$ and pick $r_1 \in Ba A_0$ so that $\|r - r_1\|_{\theta,p} <
\ep/2.$  Since $A_0$ embeds densely in $[A_0,A_1]_{\theta,p},$
there exists $m$ so that  $\|r_1 -
\frac{a_m}{\|a_{m}\|_{A_1}}\|^{\theta}_{A_1} < \ep /2 C(\theta,p)^{-1}
2^{\theta -1}.$ Using the interpolation inequality of section 2, we get that
$\|r_1 - \frac{a_m}{\|a_m\|_{A1}}\|_{\theta,p} < \ep/2.$
Therefore $R$ is separable.

Let ${\cal U}$ be the group of all isometries
of $R.$
If  $(h_n)_n$ is  a dense subset of $R,$ then $d(S,T) = \sum_{n=1}^{\infty}
2^{-n} \|Sh_n - Th_n\|_{R}$  for $S,T \in {\cal L}(R)$ is a distance that is
uniformly equivalent on ${\cal U}$ to the uniform structure of pointwise
convergence on $R,$ since ${\cal U}$    is equicontinuous. Therefore $G$ is
closed.
Finally to see that $T^{-1}$ is continuous, it suffices to repeat the same
argument used to show that G is closed.  So $X$ is homeomorphic to a group
of isometries on the reflexive space $R.$
\qed

\rmk

\noindent{If} $S$ is the symmetrization operator defined on ${\cal C}_b(X)$ by 
$Sf(x) = f(-x)$ for $x \in X$ then $S$ is an isometry  of $R$ so that
$S T_x = T_{-x}S.$

Indeed $S$ commutes with $\Phi$ and with the same argument used
for the  $T_x'$s we can define $S$ as an isometry of $R.$ Furthermore,
if $f \in R$ and $z \in X$ then $T_{-x} Sf (z) = T_{-x} f(-z) = f (-z+x) =
f(-(z-x)) = S T_x f (z).$
 
Let us now consider the space ${\cal F}_X$ of all types on a separable
stable Banach space $X.$ Denote by $\widetilde{G}$ the closure of the group
$G$ for the weak operator topology. Recall that this topology
${\tau}_{R\times R^*}$is generated by the 
family of semi-norms $p_{r,r^*} (T) = |<Tr,r^*>|$ where $r \in R$ and $r^*
\in R^*.$
Our goal here is to represent ${\cal F}_X$ in a similar way as it is  done
by S. Guerre and M. Levy in \cite{bib:G.L}. 

\begin{thm}
Let $X$ be a separable stable Banach space.
There exists a homeomorphism $T$ from $(F_X, w)$ onto $(\widetilde{G},
\tau_{R\times R^*})$ so that:
  \begin{itemize}
  \item If $\sigma_{\infty}$ is the "infinite type" then $T
(\sigma_{\infty}) = 0$
  and reciprocally, if $T \sigma = 0$ then $\sigma$ is the infinite type.
\item If $S$ is the symmetrization operator then $S \circ T_{-\sigma} =
T_{\sigma} \circ S.$
\item If $\sigma$ and $\tau$ are two types on $X$ then $T_{\sigma * \tau} =
T_{\tau * \sigma} = T_{\sigma} \circ T_{\tau}.$
\end{itemize}
\end{thm}

\pf
Let $\sigma \in F_X$ be a type defined by a sequence $(x_n)$ and an
ultrafilter ${\cal V}.$  Since $R$ is a separable reflexive space, the unit
ball of $({\cal L}(R),\tau_{R\times R^*})$ is compact metrizable, and so is 
$(\widetilde{G}, \tau_{R\times R^*}).$ Let 
$T_{\sigma}$ be the operator so that $<T_{\sigma} r , r^*> = \lim_{n,{\cal
V}} <T_{x_n} r , r^*>.$
Let $r = \Phi (\delta_{y})$ and $r^* = j^* (\delta_{z})$ where $j$ is the
natural embedding
of $[\Phi({\cal M}(X)), {\cal C}_b (X)]_{\theta, p}$ in ${\cal C}_b (X).$
The operator $T_{\sigma}$ satisfies:
\begin{eqnarray*}
<T_{\sigma} \Phi (\delta_{y}), j^* (\delta_{z})> & = & < j T_{\sigma} \Phi
(\delta_y), \delta_z> \\
    & = & \lim_{n, {\cal V}}< j T_{x_n} \Phi (\delta_y), \delta_z> \\
    & = & \lim_{n, {\cal V}} e^{- \|x_n + y - z\|} \\
    & = & e^{-\sigma (y - z)}
   \end{eqnarray*}
 \noindent{(In} case $\sigma_{\infty}$ is the infinite type,
we  have by what precedes that each $<T_{\sigma_m} \Phi (\delta_{y}), j^*
(\delta_{z})> =   e^{-\sigma_m (y - z)}$ where $\sigma_m \ra \sigma.$)

\noindent{So} to show that  the operator $T_{\sigma}$ depends only on
$\sigma$, that $T$ is continuous, and $T_{\sigma_{\infty}}  = 0$ it suffices
to show the following:

\begin{lem}
Under the same hypothesis as above,
the subset $A = \{\Phi (\delta_y)\,\, y \in X \}$ (respectively $A^*= \{
j^*(\delta_z), \,\, z \in X\}$) is total in $R$ (respectively in $R^*).$
\end{lem}

\pf
It is well known that $\{ \delta_y, \,\, y \in X \}$ is weak-* total in 
$Ba {\cal M}(X).$ By  weak compactness of $\Phi,$  the set $W = \bar{\Phi
(Ba {\cal M}(X))}$ is weakly compact in ${\cal C}_b(X)$ so $A = \{\Phi
(\delta_y)\,\, y \in X \}$ is weakly total in $W,$ which is convex, so in fact
$A = \{\Phi (\delta_y)\,\, y \in X \}$ is total  for the ${\cal C}_b(X)$ norm,
in $W.$
By density of $W$ in $Ba R,$ if $\ep > 0$ and $\psi \in BaR$ there exists
$\eta \in W$ so that $\| \psi - \eta\|_{R} < \ep/2;$ and  by totality in
norm of $A$ in $W$ there exist $(a_i)^n_{i=1}$ reals and $(y_i)_{i=1}^n$ in
$X$ so that
$\| \eta - \sum_{i=1}^n a_i \Phi(\delta_{y_i})\|_{{\cal C}_b(X)}^{\theta} < 
\ep/2 C(\theta,p)^{-1} 2^{\theta -1}.$
The totality of $A = \{\Phi (\delta_y)\,\, y \in X \}$ follows then from the
interpolation inequality of Section 2, for
$$\| \eta - \sum_{i=1}^n a_i \Phi(\delta_{y_i})\|_{\theta, p} \leq C(\theta,p)
\| \eta - \sum_{i=1}^n a_i \Phi(\delta_{y_i})\|_{\cup nW}^{1-\theta} \times
\| \eta - \sum_{i=1}^n a_i \Phi(\delta_{y_i})\|_{{\cal C}_b(X)}^{\theta}.$$
\noindent{Now} for the set $A^* =  \{ j^*(\delta_z), \,\, z \in X \},$ it is
weak-* total in $BaR^*.$ Since $R$ is reflexive, $A$ is weakly total in $BaR^*,$
which is convex so $A^*$ is in fact total in norm. This proves the lemma.
\qed

\noindent{Now} to see that $T^{-1}$ is continuous, it suffices to observe
that if $T_{\sigma_n} \ra T_{\sigma}$ for the weak operator topology then
for all $y \in X$ we have that  $$<T_{\sigma_n} \Phi (\delta_y), j^*
(\delta_0)> = e^{-\sigma_n (y)} \ra e^{-\sigma (y)}.$$

Part 2 of the Theorem follows from the identities:
\begin{eqnarray*}
<S T_{-\sigma} r, r^*> & = & \lim_{n, {\cal V}} <S T_{-x_n} r, r^*> \\
 & = & \lim_{n, {\cal V}}< T_{-x_n} r, S^* r^*> \\
 & = & \lim_{n, {\cal V}}<T_{-x_n}S^2  r,S^*r^*> \\
 & = &  \lim_{n, {\cal V}}<S T_{-x_n}S r, S^*r^*> \\
 & = &  \lim_{n, {\cal V}}< T_{-x_n}S r, S^* S^* r^*> \\  
  & = & \lim_{n, {\cal V}}< T_{-x_n}S r, r^*> \\
 & = & < T_{\sigma}  S r, r^*> 
\end{eqnarray*}

Part 3 of the Theorem follows from the following lemma and the fact that
$T_{\sigma} T_{\tau} = T_{{\sigma}*{\tau}}$ in our case.  We will prove 
the later fact after Lemma 6.
\begin{lem}
Let $G$ be an  abelian group of isometries on a reflexive Banach space $R.$ 
Let $S$ and $T$ two elements of the  closure $\widetilde{G}$ for the weak
operator topology. Then $ST = TS.$
\end{lem}

\pf 
For all $r \in R$ and $r^{*} \in  R^*$ we easily have:
\begin{eqnarray*}
<S T r,r^*> & = & <T r , S^* r^*> \\
 & = & \lim_{\alpha , {\cal U}} < T_{\alpha} r , S^* r^*> \\
 & = & \lim_{\alpha , {\cal U}} < ST_{\alpha} r ,  r^*> \\
 & = & \lim_{\alpha , {\cal U}} \lim_{\beta, {\cal V}}  <
S_{\beta}T_{\alpha} r ,  r^*> \\
 & = & \lim_{\alpha , {\cal U}} \lim_{\beta, {\cal V}}  <
T_{\alpha}S_{\beta} r ,  r^*> \\
 & = & \lim_{\beta , {\cal V}} \lim_{\alpha, {\cal U}}  <
T_{\alpha}S_{\beta} r ,  r^*> \\
 & = & \lim_{\beta , {\cal V}} \lim_{\alpha, {\cal U}}  < S_{\beta} r ,
T_{\alpha}^*r^*> \\
 & = & \lim_{\alpha , {\cal U}} \lim_{\beta, {\cal V}}  < S_{\beta} r ,
T_{\alpha}^*r^*> \\
 & = & \lim_{\alpha , {\cal U}} \lim_{\beta, {\cal V}}  <
T_{\alpha}S_{\beta} r ,  r^*> \\
& = & \lim_{\beta, {\cal V}}  < T S_{\beta} r ,  r^*> \\
& = & \lim_{\beta, {\cal V}}  <  S_{\beta} r , T^* r^*> \\
& = &  <  TS r , T^* r^*>
\end{eqnarray*}
 
This proves Lemma 6. 
To see that  $T_{\sigma} T_{\tau} = T_{\sigma * \tau}$
let  $r = \Phi (\delta_{y})$ and $r^* = j^* (\delta_{z})$ where $j$ is the
natural embedding
of $[\Phi({\cal M}(X)), {\cal C}_b (X)]_{\theta, p}$ in ${\cal C}_b (X).$ 
Let also $x_n$ and 
$y_m$ be determining sequences for the types ${\sigma}$ and  ${\tau}$
respectively.
We then have:
\begin{eqnarray*}
<  T_{\sigma} T_{\tau} r,r^*> & = & <T_{\tau}  r , (T_{\sigma} )^* r^*> \\
 & = & \lim_{m , {\cal V}} < T_{y_m} r ,  (T_{\sigma} )^*  r^*> \\
& = & \lim_{m , {\cal V}} < T_{\sigma} T_{y_m} r ,   r^*> \\
 & = & \lim_{n, {\cal U}} \lim_{m, {\cal V}}  < T_{x_n} T_{y_m} r ,   r^*> \\
& = & \lim_{n, {\cal U}} \lim_{m, {\cal V}}  < T_{x_n + y_m} r ,   r^*> \\
& = & \lim_{n, {\cal U}} \lim_{m, {\cal V}}  
< jT_{x_n + y_m}\Phi (\delta_{y}) , \delta_{z}  > \\
& = & \lim_{n, {\cal U}} \lim_{m, {\cal V}}  
\exp{-\| {x_n + y_m +y  - z}\|  }  \\
&  =  & \exp{- {\sigma} * {\tau}(y-z)} \\
&  =  & <  T_{\sigma * \tau}\Phi (\delta_{y}), j^* (\delta_{z})> 
\end{eqnarray*}
 
The fact mentioned above follows then from Lemma 5.
\qed

Let us now turn to  comparing the strong topology on types and the strong
operator topology $\tau_s$ on the representation group $G.$  Let us recall
that Y. Raynaud proved in \cite{bib:R} that the topology $s$  of uniform
convergence on bounded subsets of $X$ on  the space of types 
${\cal F}_X$ is separable  provided the space $X$ is stable separable.
 
\begin{thm}
Let $X$ be a separable stable Banach space and let $T$ be the representation 
of $F_X.$ If $(\sigma_{n})$ is a sequence of types in $F_X$ then
  $(\sigma_{n}$ converges for the topology $s$  to    a type $\sigma$ 
if and only if $T_{\sigma_{n}}$ converges strongly to $T_{\sigma}$.
 \end{thm}

\pf
Suppose that $\sigma_n \ra \sigma$ for the topology $s.$ This means that 
for any positive real $M$ we have that $d_M (\sigma_n, \sigma) = \sup_{\|x\|
\leq M}|\sigma_n(x) - \sigma (x) | \ra 0$ as $n \ra \infty.$
Consider now the contractions  $T_{\sigma_{n}}$ and $T_{\sigma}.$
Since we have seen before that the family $\Phi (\delta_{x})$
was total in $X,$ it suffices to show that $T_{\sigma_{n}}\Phi (\delta_{x})=
T_{\sigma_{n}} T_x \Phi (\delta_0) = T_x T_{\sigma_{n}} (\Phi (\delta_0))$
converges in $R,$ i.e. if $U_x = \Phi (\delta_{x})$ then $T_{\sigma_{n}}
U_0$ converges to $T_{\sigma} U_0$ in $R.$
Now each $T_{\sigma_{n}}$ and $T_{\sigma}$ are weak limits of elements in
$\Phi({\cal M}(X)),$ which is convex, so its weak closure is the same as its
closure for the norm defined by its gauge. Therefore 
$\| T_{\sigma_{n}} U_0 - T_{\sigma} U_0\|_{A_0} \leq 2.$
On the other hand,
\begin{eqnarray*}
\| T_{\sigma_{n}} U_0 - T_{\sigma} U_0\|_{{\cal C}_b (X)} & = & \sup_{ y \in
X}|( T_{\sigma_{n}} U_0 - T_{\sigma} U_0) (y)| \\
& = & \sup_{ y \in X}|< j(T_{\sigma_{n}} U_0 - T_{\sigma} U_0), \delta_y >| \\
& = & \sup_{ y \in X}|< T_{\sigma_{n}} U_0 - T_{\sigma} U_0, j^*(\delta_y) >| \\
& = & \sup_{ y \in X}| e^{-\sigma_n (-y)} - e^{-\sigma (-y)}|
\end{eqnarray*}
\noindent{But} the last quantity tends to zero as $n$ tends to infinity since
on the bounded sets of $X$ that is the convergence in the $s$ sense, while
outside the bounded sets, this is a property of the exponential function.
We then conclude the proof of one direction using the interpolation inequality
$$\| T_{\sigma_{n}} U_0 - T_{\sigma} U_0\|_{R} \leq C(\theta, p) 2^{1 - \theta} 
\| T_{\sigma_{n}} U_0 - T_{\sigma} U_0\|_{{\cal C}_b (X)}^{\theta}.$$
For the converse, if $T_{\sigma_n} \ra T_{\sigma}$ in $R$  as $n \ra \infty,$
then $j(T_{\sigma_n} U_0)   \ra j(T_{\sigma} U_0)$ in ${\cal C}_b(X).$ 
But for that norm we clearly have:
\begin{eqnarray*}
\|j(T_{\sigma_n} U_0) \ra j(T_{\sigma} U_0)\|_{{\cal C}_b(X)} & \geq &
|<(T_{\sigma_{n}} U_0 - T_{\sigma} U_0) , \delta_{-x}>| \\
& = & |e^{-\sigma_{n}(x)} - e^{-\sigma(x)}|\\
& = & |e^{-\sigma (x)} [e^{-\sigma_n (x) + \sigma(x)} - 1]|\\
& \geq & e^{-M - \|\sigma\|} |e^{-\sigma_n (x) + \sigma(x)} - 1| 
\end{eqnarray*}
\noindent{if} we suppose that $\|x\| \leq M.$
Thus $d_M(\sigma_n, \sigma) = \sup_{\|x\| \leq M} |\sigma_n (x) - \sigma (x)
| \ra 0$ as $n \ra \infty.$
\qed

\bigskip

School of Sciences and Engineering

Al Akhawayn University in Ifrane

Ifrane 53000 Morocco.

\end{document}